\documentclass[final]{siamltex}
\usepackage{epsfig}
\usepackage{graphicx}
\usepackage{amsmath}
\usepackage{amssymb}
\usepackage{setspace}
\usepackage{graphicx}
\usepackage{graphicx}

\newtheorem{defn}{\noindent $\mathbf{Definition}$}[section]

\newtheorem{thm}[defn]{$\mathbf{Theorem}$}

\title{Convergence of an Iterative Algorithm for Teichm\"uller Maps via Generalized Harmonic Maps}

\author{Lok Ming Lui, Xianfeng Gu and Shing-Tung Yau}

\begin{document}

\maketitle

\begin{abstract}
Finding surface mappings with least distortion arises from many applications in various fields. Extremal Teichm\"uller maps are surface mappings with least conformality distortion. The existence and uniqueness of the extremal Teichm\"uller map between Riemann surfaces of finite type are theoretically guaranteed \cite{Alastair}. Recently, a simple iterative algorithm for computing the Teichm\"uller maps between connected Riemann surfaces with given boundary value was proposed in \cite{LuiTMap}. Numerical results was reported in the paper to show the effectiveness of the algorithm. The method was successfully applied to landmark-matching registration. The purpose of this paper is to prove the iterative algorithm proposed in \cite{LuiTMap} indeed converges.
\end{abstract}



\pagestyle{myheadings}
\thispagestyle{plain}
\markboth{Lui, Gu and Yau}{Teichm\"uller Maps via Harmonic Energy Optimization}

\section{Introduction}
Finding meaningful surface mappings with least distortion has fundamental importance. Applications can be found in different areas such as registration, shape analysis and grid generation. Conformal mapping has been widely used to establish a good one-to-one correspondence between different surfaces, since it preserves the local geometry well \cite{Haker,Fischl2, Gu1, Gu2,Gu3,Hurdal,conformal1,conformal2,conformal3}. The Riemann mapping theorem guarantees the existence of conformal mappings between simply-connected surfaces. However, this fact is not valid for general Riemann surfaces. Given two Riemann surfaces with different conformal modules, there is generally no conformal mapping between them. In this case, it is usually desirable to obtain a mapping that minimizes the conformality distortion. Every diffeomorphic surface mapping is associated with a unique Beltrami differential, which is a complex-valued function, $\mu_f$, defined on the source surface. The Beltrami differential, $\mu_f$, measures the deviation of the mapping from a conformal map. Given two Riemann surfaces $S_1$ and $S_2$, there exists a unique and bijective map $f:S_1\to S_2$, called the {\it Teichm\"uller map}, minimizing the $L^{\infty}$ norm of the Beltrami differential \cite{Alastair}. Therefore, the extremal Teichm\"uller map can be considered as the `most conformal' map between Riemann surfaces of the same topology, which is a natural extension of conformal mappings.

\subsection{Extremal problem}
Mathematically, the extremal problem for obtaining a surface mapping with least conformality distortion can be formulated as follows. Suppose $(S_1, \sigma |dz|^2)$ and $(S_2, \rho |dw|^2)$ are two Riemann surfaces of finite type, where $z$ and $w$ are their conformal parameters respectively. Every diffeomorphism between $S_1$ and $S_2$ is associated with a unique Beltrami differential. A Beltrami differential $\mu(z) \frac{d\bar{z}}{dz}$ on the Riemann surface $S_1$ is an assignment to each chart $(U_{\alpha},\phi_{\alpha})$ of an $L^{\infty}$ complex-valued function $\mu_{\alpha}$, defined on local parameter $z_{\alpha}$. Then, $f:S_1 \to S_2$ is said to be a quasi-conformal mapping associated with the Beltrami differential $\mu (z) \frac{\overline{dz}}{dz}$ if for any chart $(U_{\alpha},\phi_{\alpha})$ on $S_1$ and any chart $(V_{\beta},\psi_{\beta})$ on $S_2$, the mapping $f_{\alpha \beta}:= \psi_{\beta}\circ f\circ {\phi}_{\alpha}^{-1}$ is quasi-conformal associated with $\mu_{\alpha}(z_{\alpha})\frac{d\overline{z_{\alpha}}}{dz_{\alpha}}$.

Our goal is to look for an extremal quasi-conformal mapping, which are extremal in the sense of minimizing the $||\cdot||_{\infty}$ over all Beltrami differentials corresponding to quasi-conformal mappings between $S_1$ and $S_2$. The idea of extremality is to make the supreme norm of the Beltrami differential as small as possible such that $f$ is as `nearly conformal' as possible.

The extremal problem can therefore be formulated as finding $f:S_1 \to S_2$ that solves:
\begin{equation}\label{extremal}
f = \mathbf{argmin}_{g\in \mathcal{A}} \{||\mu_g||_{\infty}\}
\end{equation}

\noindent where $\mathcal{A}=\{g:S_1 \to S_2: g \mathrm{\ is\ a\ diffeomorphism}\}$.

The above optimization problem (\ref{extremal}) has a unique global minimizer provided that $S_1$ and $S_2$ are Riemann surfaces of finite type. Also, the unique minimizer $f: S_1 \to S_2$ is a Teichm\"uller map, that is, its associated Beltrami differential $\mu_{f}$ is of the following form:
\begin{equation}
\mu_{f} = k \frac{\bar{\varphi}}{|\varphi|}
\end{equation}
\noindent where $0 \leq k <1$ is a non-negative real constant and $\varphi$ is an integrable holomorphic function defined on $S_1$ ($\varphi \neq 0$). Beltrami differential of this form is said to be of Teichm\"uller type.

\subsection{An iterative algorithm for Teichm\"uller maps}
To solve the extremal problem (\ref{extremal}) to obtain the Teichm\"uller map between connected surfaces, an iterative algorithm was proposed in \cite{LuiTMap}, called the quasi-conformal(QC) iteration. The ultimate goal is to obtain the extremal map between connected (either simply-connected or multiply-connected) surfaces with given boundary value, which minimizes the conformality distortion. The basic idea of the iterative algorithm is to project the Beltrami differential to the space of all Beltrami differentials of Teichm\"uller type, and compute a quasi-conformal map whose Beltrami differential is closest to the projection in the least square sense. More specifically, the QC iteration for solving (\ref{extremal}) can be described as follows.
\begin{equation}\label{QCiteration}
\begin{cases}
f_{n+1} = \mathbf{LBS}(\mu_{n+1}),\\
\widetilde{\mu}_{n+1} = \mu_n + \alpha \mu(f_{n+1},\mu_n),\\
\mu_{n+1} =\mathcal{L}(\mathcal{P}(\widetilde{\mu}_{n+1}))\\
\end{cases}
\end{equation}
\noindent where $f_{n}$ is the quasi-conformal map obtained at the $n^{th}$ iteration, $\nu_n$ is the Beltrami differential of $f_n$ and $\mu_n$ is a Beltrami differential of constant modulus.

$\mathbf{LBS}(\mu)$ is the operator to obtain a quasi-conformal mapping whose Beltrami differential is closest to $\mu$ in the least square sense. In other words,
\begin{equation}
\mathbf{LBS}(\mu) = \mathbf{argmin}_{f\in \mathcal{A}} \{ \int_{S_1} |\frac{\partial f}{\partial \bar{z}} - \mu \frac{\partial f}{\partial z}|^2 dS_1\}
\end{equation}

$\mu(f_{n+1},\nu_n)$ denotes the Beltrami differential of $f_{n+1}$ under the auxiliary metric with respect to $\nu_n$, namely, $|dz + \nu_n d\bar{z}|^2$ ($|dz|^2$ is the original metric on $S_1$). More precisely, $\mu(f_{n+1},\nu_n)$ can be explicitly computed as follows:
\begin{equation}
\mu(f_{n+1},\nu_n) = \left(\frac{\partial f_{n+1}}{\partial \bar{z}} + \nu_{n} \frac{\partial f_{n+1}}{\partial z}\right)/\left(\frac{\partial f_{n+1}}{\partial z} - \overline{\nu_{n}} \frac{\partial f_{n+1}}{\partial \bar{z}}\right)
\end{equation}

$\mathcal{P}(\widetilde{\mu}_{n+1})$ is the operator to project $\widetilde{\mu}_{n+1}$ to the space of Beltrami differentials with constant modulus. It is defined as:
\begin{equation}
\mathcal{P}(\widetilde{\mu}_{n+1}) = \mu_n + \epsilon w_n
\end{equation}
\noindent where $w_n:S_1\to \mathbb{C}$ and $\epsilon: S_1\to \mathbb{R}^+$ is a suitable real function on $S_1$ such that $|\mu_n + \epsilon w_n|$ is a constant.

In practice, the projection operator can be simplified as
\begin{equation}
\mathcal{P}(\widetilde{\mu}_{n+1}) = \left( \frac{\int_{S_1} |\widetilde{\mu}_{n+1} | dS_1}{\int_{S_1} dS_1}\right) \frac{\widetilde{\mu}_{n+1} }{|\widetilde{\mu}_{n+1} |}
\end{equation}

\noindent $\mathcal{L}$ is the Laplacian smoothing operator to smooth out $\mathcal{P}(\widetilde{\mu}_{n+1} )$.

Both $\mathcal{P}(\nu_n)$, $\mathbf{LBS}(\mu_{n+1})$ and  $\mu(f_{n+1},\nu_n)$ can be easily computed. In particular, the discretization of $\mathbf{LBS}(\mu_{n+1})$ on a triangulation mesh can be reduced to a least square problem of a linear system.

When $\mu(f_{n+1},\mu_n)$ is small, $\alpha$ can be chosen to be 1. Then, $\mu_n + \alpha \mu(f_{n+1},\mu_n) \approx \mu(f_{n+1},0)$, where $\mu(f_{n+1},0)$ is the Beltrami differential of $f_{n+1}$ under the original metric.. The QC iteration can be further modified as
\begin{equation}\label{QCiteration2}
\begin{cases}
f_{n+1} = \mathbf{LBS}(\mu_{n}),\\
\widetilde{\mu}_{n+1} = \mu(f_{n+1},0),\\
\mu_{n+1} =\mathcal{P}(\mathcal{L}(\widetilde{\mu}_{n+1})).
\end{cases}
\end{equation}

The QC iteration (\ref{QCiteration}) is very efficient. Also, numerical results reported in \cite{LuiTMap} demonstrate that the proposed iteration can compute the Teichm\"uller map accurately, even on highly irregular meshes. The algorithm was successfully applied to landmark-based registration for applications in medical imaging and computer graphics.

This paper is to provide a complete analysis of the above iterative algorithm (\ref{QCiteration}). In particular, we prove the convergence of (\ref{QCiteration}) that $f_n$ and $\mu_n$ respectively converges to the extremal map $f^*$ and its associated Beltrami differential $\mu^*$, which solves the optimization problem (\ref{extremal}).

We remark that although the iterative algorithm is designed for obtaining extremal map between connected surfaces with given boundary values, the convergence proof applies to general Riemann surfaces of finite type (such as high-genus closed surfaces). In other words, the QC iteration can be applied to computing extremal maps between general Riemann surfaces. For the ease of the presentation, we will restrict our discussion to the situation when both $S_1$ and $S_2$ are either simply-connected or multiply-connected open surfaces.

\subsection{Organization}
This paper is organized as follows. In Section \ref{sec:background}, we describe some mathematical background, which is relevant to this work. In Section \ref{sec:QCiteration}, we reformulate the extremal problem defined by (\ref{extremal}) as the optimization problem of the harmonic energy, which helps us to understand the iterative algorithm (\ref{QCiteration}) better. In Section \ref{sec:convergence}, we prove the convergence of the QC iteration to our desired extremal Teichm\"uller map. A concluding remark will be given in Section \ref{sec:conclusion}.

\section{Mathematical background}\label{sec:background}
\subsection{Quasi-conformal mappings and Beltrami equation}
In this section, we describe some basic mathematical concepts relevant to our algorithms. For details, we refer the readers to \cite{Gardiner,Lehto}.

A surface $S$ with a conformal structure is called a \emph{Riemann surface}. Given two Riemann surfaces $S_1$ and $S_2$, a map $f:S_1\to S_2$ is \emph{conformal} if it preserves the surface metric up to a multiplicative factor called the {\it conformal factor}. A generalization of conformal maps is the \emph{quasi-conformal} maps, which are orientation preserving homeomorphisms between Riemann surfaces with bounded conformality distortion, in the sense that their first order approximations takes small circles to small ellipses of bounded eccentricity \cite{Gardiner}. Mathematically, $f \colon \mathbb{C} \to \mathbb{C}$ is quasi-conformal provided that it satisfies the Beltrami equation:
\begin{equation}\label{beltramieqt}
\frac{\partial f}{\partial \overline{z}} = \mu(z) \frac{\partial f}{\partial z}.
\end{equation}
\noindent for some complex-valued function $\mu$ satisfying $||\mu||_{\infty}< 1$. $\mu$ is called the \emph{Beltrami coefficient}, which is a measure of non-conformality. $\mu_f$ measures how far the map is deviated from a conformal map. $\mu \equiv 0$ if and only if $f$ is conformal. Infinitesimally, around a point $p$, $f$ may be expressed with respect to its local parameter as follows:
\begin{equation}
\begin{split}
f(z) & = f(p) + f_{z}(p)z + f_{\overline{z}}(p)\overline{z} \\
& = f(p) + f_{z}(p)(z + \mu(p)\overline{z}).
\end{split}
\end{equation}

\begin{figure*}[t]
\centering
\includegraphics[height=1.5in]{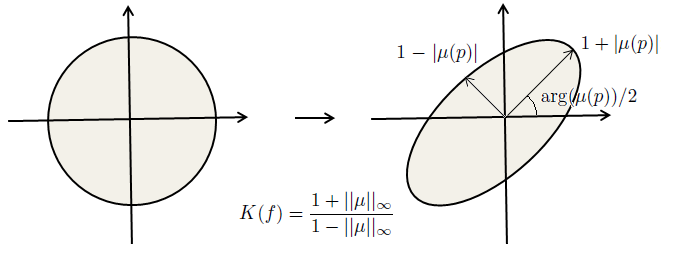}
\caption{Illustration of how the Beltrami coefficient determines the conformality distortion. \label{fig:illustration1}}
\end{figure*}
Obviously, $f$ is not conformal if and only if $\mu(p)\neq 0$. Inside the local parameter domain, $f$ may be considered as a map composed of a translation to $f(p)$ together with a stretch map $S(z)=z + \mu(p)\overline{z}$, which is postcomposed by a multiplication of $f_z(p)$. All the conformal distortion of $S(z)$ is caused by $\mu(p)$. $S(z)$ is the map that causes $f$ to map a small circle to a small ellipse. From $\mu(p)$, we can determine the directions of maximal magnification and shrinking and the amount of their distortions as well. Specifically, the angle of maximal magnification is $\arg(\mu(p))/2$ with magnifying factor $1+|\mu(p)|$; The angle of maximal shrinking is the orthogonal angle $(\arg(\mu(p)) -\pi)/2$ with shrinking factor $1-|\mu(p)|$. Thus, the Beltrami coefficient $\mu$ gives us all the information about the properties of the map (see Figure \ref{fig:illustration1}).

The maximal dilation of $f$ is given by:
\begin{equation}\label{dilation}
K(f) = \frac{1+||\mu||_{\infty}}{1-||\mu||_{\infty}}.
\end{equation}

Quasiconformal mapping between two Riemann surfaces $S_1$ and $S_2$ can also be defined. Instead of the Beltrami coefficient, the {\it Beltrami differential} has to be used. A Beltrami differential $\mu(z) \frac{d\bar{z}}{dz}$ on the Riemann surface $S_1$ is an assignment to each chart $(U_{\alpha},\phi_{\alpha})$ of an $L^{\infty}$ complex-valued function $\mu_{\alpha}$, defined on local parameter $z_{\alpha}$ such that
\begin{equation}
\mu_{\alpha}(z_{\alpha})\frac{d\overline{z_{\alpha}}}{dz_{\alpha}} = \mu_{\beta}(z_{\beta})\frac{d\overline{z_{\beta}}}{dz_{\beta}},
\end{equation}
\noindent on the domain which is also covered by another chart $(U_{\beta},\phi_{\beta})$, where $\frac{dz_{\beta}}{dz_{\alpha}}= \frac{d}{dz_{\alpha}}\phi_{\alpha \beta}$ and $\phi_{\alpha \beta} = \phi_{\beta}\circ \phi_{\alpha}^{-1}$.

\begin{figure*}[t]
\centering
\includegraphics[height=1.85in]{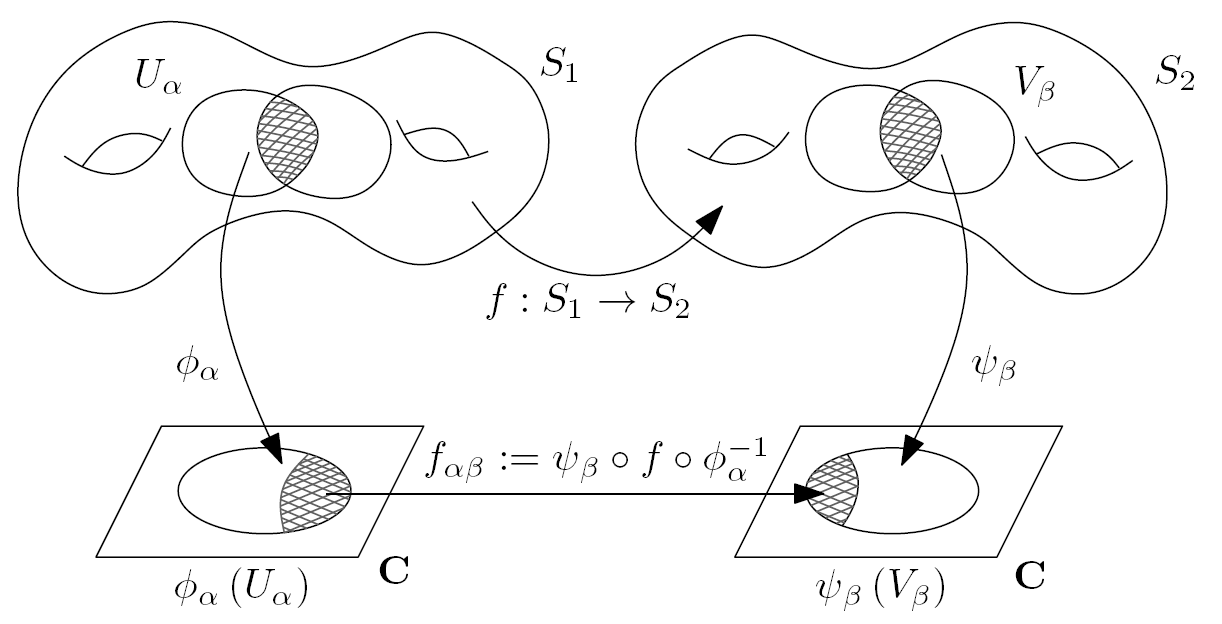}
\caption{Illustration of the definition of quasi-conformal map between Riemann surfaces. \label{fig:qcillustration2}}
\end{figure*}

An orientation preserving diffeomorphism $f:S_1 \to S_2$ is called quasi-conformal associated with $\mu(z) \frac{\overline{dz}}{dz}$ if for any chart $(U_{\alpha},\phi_{\alpha})$ on $S_1$ and any chart $(V_{\beta},\psi_{\beta})$ on $S_2$, the mapping $f_{\alpha \beta}:= \psi_{\beta}\circ f\circ {\phi}_{\alpha}^{-1}$ is quasi-conformal associated with $\mu_{\alpha}(z_{\alpha})\frac{d\overline{z_{\alpha}}}{dz_{\alpha}}$. See Figure \ref{fig:qcillustration2} for an illustration.

\subsection{Extremal maps and Teichm\"uller maps}
A special class of quasi-conformal maps is called the {\it extremal} maps, which minimize the conformality distortion.  More specifically, an extremal quasi-conformal map between $S_1$ and $S_2$ is extremal in the sense of minimizing the $||\cdot||_{\infty}$ over all Beltrami differentials corresponding to quasi-conformal mappings between the two surfaces. Extremal map always exists but need not to be unique. Mathematically, an extremal quasi-conformal mapping can be defined as follows:

\bigskip

\begin{defn}
Suppose $S_1$ and $S_2$ are connected Riemann surfaces with boundaries. Let $f: S_1\to S_2$ be a quasi-conformal mapping between $S_1$ and $S_2$. $f$ is said to be an {\it extremal map} if for any quasi-conformal mapping $h: S_1\to S_2$ isotopic to $f$ relative to the boundary,
\begin{equation}\label{inequality}
K(f) \leq K(h)
\end{equation}
\noindent It is uniquely extremal if the inequality (\ref{inequality}) is strict.
\end{defn}

Closely related to the extremal map is the {\it Teichm\"uller} map. It is defined as follows.

\bigskip

\begin{defn}
Let $f:S_1\to S_2$ be a quasi-conformal mapping. $f$ is said to be a Teichm\"uller map associated to the integrable holomorphic function $\varphi:S_1\to \mathbb{C}$ if its associated Beltrami differential is of the form:
\begin{equation}\label{Teichmullermap}
\mu(f) = k \frac{\overline{\varphi}}{|\varphi|}
\end{equation}
\noindent for some constant $k <1$ and holomorphic function $\varphi \neq 0$ with $||\varphi||_1 = \int_{S_1} |\varphi| <\infty$.
\end{defn}

In other words, a Teichm\"uller map is a quasi-conformal mapping with uniform conformality distortion over the whole domain.

\bigskip

Extremal map might not be unique. However, a Teichm\"uller map associated with a holomorphic function is the unique extremal map in its homotopic class. In particular, a Teichmuller map between two connected open surfaces with suitable given boundary values is the unique extremal map. The Strebel's theorem explains the relationship bewtween the Teichm\"uller map and extremal map.

\bigskip

\begin{defn}[Boundary dilation]
The {\it boundary dilation} $K_1[f]$ of $f$ is defined as:
\begin{equation}
K_1[f] = \inf_{C} \{K(h|_{S_1 \setminus C}): h \in \mathfrak{F}, C \subseteq S_1, C \mathrm{\ is\ compact.}\}
\end{equation}
\noindent where $\mathfrak{F}$ is the family of quasi-conformal homeomorphisms of $S_1$ onto $S_2$ which are homotopic to $f$ modulo the boundary.
\end{defn}

\bigskip

\begin{thm}[Strebel's theorem, See \cite{Strebel}, page 319]\label{Strebel}
Let $f$ be an extremal quasi-conformal map with $K(f) >1$. If $K_1[f]< K(f)$, then $f$ is a Teichm\"uller map associated with an integrable holomorphic function on $S_1$. Hence, $f$ is also an unique extremal mapping.
\end{thm}

\bigskip

In other words, an extremal map between $S_1$ and $S_2$ with suitable boundary condition is a Teichm\"uller map. In particular, the Teichm\"uller mapping and extremal mapping of the unit disk are closely related.

\bigskip

\begin{thm}[See \cite{Reich}, page 110]\label{diskteichmuller}
Let $g:\partial \mathbb{D} \to \partial \mathbb{D}$ be an orientation-preserving homeomorphism of $\partial \mathbb{D}$. Suppose further that $h'(e^{i\theta}) \neq 0$ and $h''(e^{i\theta})$ is bounded. Then there is a Teichm\"uller map $f$ that is the unique extremal extension of $g$ to $\mathbb{D}$. That is, $f:\mathbb{D}\to \mathbb{D}$ is an extremal mapping with $f|_{\partial \mathbb{D}} = g$.
\end{thm}

\bigskip

Thus, if the boundary correspondence satisfies certain conditions on its derivatives, the extremal map of the unit disk must be a Teichm\"uller map.

Now, in the case when interior landmark constraints are further enforced, the existence of unique Teichm\"uller map can be guaranteed if the boundary and landmark correspondence satisfy suitable conditions. The unique Teichm\"uller map is extremal, which minimizes the maximal conformality distortion. The following theorem can be derived immediately from the Strebel's Theorem (Theorem \ref{Strebel}):

\bigskip

\begin{thm}\label{landmarkteichmuller}
Let $\{p_i\}_{i=1}^n \in S_1$ and $\{q_i\}_{i=1}^n \in S_2$ be the corresponding interior landmark constraints. Let $f:S_1\setminus \{p_i\}_{i=1}^n \to S_2\setminus \{q_i\}_{i=1}^n$ be the extremal map, such that $p_i$ corresponds to $q_i$ for all $1 \leq i\leq n$. If $K_1[f]< K(f)$, then $f$ is a Teichm\"uller map associated with an integrable holomorphic function on $S_1\setminus \{p_i\}_{i=1}^n$. Hence, $f$ is an unique extremal map.
\end{thm}

In particular, a unique Teichm\"uller map $f:\mathbb{D}\to \mathbb{D}$ between unit disks with interior landmark constraints enforced exists, if the boundary map $f|_{\partial \mathbb{D}}$ satisfies suitable conditions. The following theorem can be obtained directly from Theorem \ref{diskteichmuller}:

\bigskip

\begin{thm}\label{landmarkteichmullerdisk}
Let $g:\partial \mathbb{D} \to \partial \mathbb{D}$ be an orientation-preserving homeomorphism of $\partial \mathbb{D}$. Suppose further that $h'(e^{i\theta}) \neq 0$ and $h''(e^{i\theta})$ is bounded. Let $\{p_i\}_{i=1}^n \in \mathbb{D}$ and $\{q_i\}_{i=1}^n \in \mathbb{D}$ be the corresponding interior landmark constraints. Then there is a Teichm\"uller map $f:\mathbb{D}\setminus \{p_i\}_{i=1}^n \to \mathbb{D}\setminus \{q_i\}_{i=1}^n$ matching the interior landmarks, which is the unique extremal extension of $g$ to $\mathbb{D}$. That is, $f:\mathbb{D}\setminus \{p_i\}_{i=1}^n \to \mathbb{D}\setminus \{q_i\}_{i=1}^n$ is an extremal Teichm\"uller map with $f|_{\partial \mathbb{D}} = g$ matching the interior landmarks.
\end{thm}

\subsection{Harmonic maps}
Our iterative algorithm to compute Teichm\"uller maps is closely related to harmonic maps. Let $(S_1,\sigma|dz|^2)$ and $(S_2,\rho|dw|^2)$ be two Riemann surfaces of finite type, where $z$ and $w$ refer to the local conformal coordinate on the surface $S_1$ and $S_2$.

For a Lipschitz map $f:(S_1,\sigma|dz|^2)\to (S_2,\rho|dw|^2)$, we
define the energy $E(f;\sigma, \rho)$ of the map $w$ to be
\begin{equation}
    E_{harm}(f;\sigma,\rho)=\int_{S_1} \frac{1}{2} \|df\|^2 dv(\sigma) = \int_{S_1}
    \frac{\rho(w(z))}{\sigma(z)} (|w_z|^2+|w_{\bar{z}}|^2)
    \sigma(z)dzd\bar{z}.
\end{equation}
Therefore
\begin{equation}
E_{harm}(f;\sigma,\rho)=\int_{S_1}
    \rho(w(z)) (|w_z|^2+|w_{\bar{z}}|^2)dzd\bar{z}.
\end{equation}
It depends on the metric structure of the target surface $\rho|dw|^2$ and
the conformal structure $\sigma|dz|^2$ of the source.

A critical point of this functional is called a harmonic map. We
will focus on the situation where we have fixed the homotopy class
$f_0:S_1\to S_2$ of maps into the compact target $S_2$ with non-positive
curvature $K(w)\le 0$ everywhere. In that case, there is a unique
harmonic map $f(\sigma,\rho):(S_1,\sigma)\to (S_2,\rho)$ in the homotopy
class of $f_0$. If $f$ is harmonic, then
\begin{equation}
    f_{z\bar{z}} + (\log \rho)_z f_zf_{\bar{z}} \equiv 0.
\end{equation}

The pull back metric on $S_1$ induced by $f$ is given by
\begin{equation}
    f^*(\rho(w) |dw|^2) = \rho ( f_z dz + f_{\bar{z}} d\bar{z} ) ( \bar{f_z} d\bar{z} + \overline{f_{\bar{z}}} dz )
\end{equation}
Then the Hopf differential is
\begin{equation}
    \Phi(f):= \rho(f(z)) f_z \overline{f_{\bar{z}}} dz^2.
\end{equation}
It can be shown that $f$ is harmonic if and only if its Hopf differential is a holomorphic quadratic differential.

\section{Quasi-conformal iteration}\label{sec:QCiteration}
Before giving a complete analysis of the convergence of the QC iteration, we reformulate the extremal problem (\ref{extremal}) as the optimization problem of the harmonic energy, in order to better understand the iterative algorithm.

Consider two connected open surfaces $S_1$ and $S_2$ with boundaries, which are of the same topology. $S_1$ and $S_2$ can either be simply-connected or multiply-connected. Suppose $\sigma |dz|^2$ and $\rho |dw|^2$ are the Riemannian metric on $S_1$ and $S_2$ respectively. Assume $(S_2, \rho |dw|^2$) has non-positive Gaussian curvature $K(w)$ everywhere. Let $f:S_1 \to S_2$ be any quasi-conformal mapping between $S_1$ and $S_2$. In the homotopic class $[f]$ of $f$, there exists a unique Teichm\"uller map, $f^*$. $f^*$ is also extremal within the homotopic class $[f]$. More specifically, the homotopic class $[f]$ can be defined as:
\begin{equation}
[f]= \{g:S_1:S_2 : g|_{\partial S_1} = f|_{\partial S_2}\}.
\end{equation}

We have, $||\mu_{f^*}||_{\infty}\leq ||\mu_{g}||_{\infty}$ for all $g \in [f]$, where $\mu_{f^*}$ and $\mu_g$ are the Beltrami differentials of $f^*$ and $g$ respectively.

Consider the space of all admissible Beltrami differentials on $S_1$, which is denoted by $\mathcal{B}(S_1,S_2)$. Every Beltrami differential $\mu \in \mathcal{B}(S_1,S_2)$ induces a conformal structure $g(\mu)$ on $S_1$, namely,
\begin{equation}
g(\mu) = |dz + \mu d\bar{z}|^2
\end{equation}
Suppose $\mu_1, \mu_2 \in \mathcal{B}(S_1,S_2)$, we say that they are {\it globally equivalent}, if there is a biholomorphic mapping $f:(S_1,g(\mu_1)) \to (S_1, g(\mu_2))$ such that $f$ is homotopic to the identity map of $S_1$. The equivalence class of $\mu$ is represented by $[\mu]$. Each global equivalence class of Beltrami differentials has a unique representative of Teichm\"uller form. We denote the space of all Beltrami differentials of Teichm\"uller form by
\begin{equation}
\mathcal{T}(S_1,S_2) := \{\mu \in \mathcal{B}(S_1,S_2): |\mu| \mathrm{\ is\ a\ constant.}\}.
\end{equation}

We can now define an energy functional $E_{BC}$ on $\mathcal{B}(S_1,S_2)$. For any $\mu \in \mathcal{B}(S_1,S_2)$, there exists a unique harmonic map $f(\mu,\rho): (S_1, g(\mu))\to (S_2, \rho |dw|^2)\in [f]$ , which is solely determined by $\mu$ and $\rho |dw|^2$. The value of $E_{BC}(\mu)$ can then be defined as the harmonic energy of $f(\mu,\rho)$. That is,
\begin{equation}
E_{BC}(\mu) = E_{harm}(f(\mu,\rho)) = \int_{S_1} \frac{1}{2}||df(\mu,\rho)||^2
\end{equation}

\noindent $E_{BC}: \mathcal{B}(S_1,S_2)\to \mathbb{R}$ is a smooth function.

\bigskip

\begin{lemma}\label{Ebound}
The energy functional $E_{BC}:\mathcal{T}(S_1,S_2) \to \mathbb{R}$ is bounded below by
\begin{equation}
E_{BC}(\mu) \geq \int_{S_2} \rho(w) du dv
\end{equation}
\noindent where $w = u+iv$. The equality holds if and only if $(S_1,g(\mu))$ is conformally equivalent to $(S_2,\rho |dw|^2)$. And the harmonic map $f(\mu,\rho): (S_1, g(\mu))\to (S_2, \rho |dw|^2)$ is a conformal mapping.
\end{lemma}
\begin{proof}
Let $z = x+ i y$ be the local coordinate of $(S_1, g(\mu))$. The Jacobian of the mapping $f(\mu,\rho):(S_1, g(\mu))\to (S_2, \rho |dw|^2)$ is given by
\begin{equation}
J(z) = |w_z|^2 - |w_{\bar{z}}|^2.
\end{equation}
Therefore,
\begin{equation}
J(z)dxdy = (|w_z|^2 - |w_{\bar{z}}|^2) dxdy = dudv.
\end{equation}
The harmonic energy is given by
\begin{equation}
\begin{split}
E_{harm}(f(\mu,\rho)) = E_{BC}(\mu) &= \int_{S_1} \rho(w) (|w_z|^2 + |w_{\bar{z}}|^2) dxdy \\
&= \int_{S_2} \rho(w) \frac{|w_z|^2 + |w_{\bar{z}}|^2}{|w_z|^2 - |w_{\bar{z}}|^2} dudv,
\end{split}
\end{equation}
\noindent where
\begin{equation}
\begin{split}
\frac{|w_z|^2 + |w_{\bar{z}}|^2}{|w_z|^2 - |w_{\bar{z}}|^2} &= \frac{1 + |\frac{w_{\bar{z}}}{w_z}|^2}{1 - |\frac{w_{\bar{z}}}{w_z}|^2} = \frac{1+|\mu|^2}{1-|\mu|^2}\\
&= \frac{1+k^2}{1-k^2} = \frac{1}{2}\left(\frac{1+k}{1-k} + \frac{1-k}{1+k}\right)= \frac{1}{2}\left(K + \frac{1}{K}\right)\mathrm{\ \ and}
\end{split}
\end{equation}
\begin{equation}
k = |\mu|, \ 0 \leq k \leq 1,\  K = \frac{1+k}{1-k},\ K\geq 1.
\end{equation}
Hence,
\begin{equation}
E_{BC}(\mu) = \frac{1}{2}\int_{S_2} \rho(w) \left(K + \frac{1}{K}\right) dudv \geq \frac{1}{2}\int_{S_2} \rho(w) (2) dudv  = \int_{S_2}\rho(w) dudv.
\end{equation}
Equality holds if and only if $K \equiv 1$, namely, $k\equiv 0$. This implies $f(\mu,\rho)$ is a conformal mapping.
\end{proof}

\bigskip

\begin{theorem}
The global minimizer of the energy functional $E_{BC}:\mathcal{T}(S_1,S_2) \to \mathbb{R}$ is the Beltrami differerntial associated to the unique Teichm\"uller map between $(S_1, \sigma |dz|^2)$ and $(S_2, \rho |dw|^2)$ in the homotopic class $[f]$ of $f$.
\end{theorem}
\begin{proof}
Let $\mu^*$ be the Beltrami differential of the Teichm\"uller map $\tilde{f}$. It suffices to show that $\tilde{f} :(S_1, g(\mu^*))\to (S_2, \rho |dw|^2)$ is a conformal mapping.

To see this, let $\tilde{f}^*(\rho |dw|^2)$ denote the pull back metric. Then,
\begin{equation}
    \tilde{f}^*(\rho |dw|^2)= e^{2\lambda_2( \tilde{f}(z) )} |df(z)|^2.
\end{equation}

Under the pull back metric, the mapping $\tilde{f}:( S_1, \tilde{f}^*(\rho |dw|^2)) \to (S_2,\rho |dw|^2)$ is isometric. We have
\begin{equation}
\begin{array}{lcl}
    d\tilde{f}(z) &=& \frac{\partial \tilde{f}(z)}{\partial z} dz + \frac{\partial \tilde{f}(z)}{\partial \bar{z}} d \bar{z}\\
    &=& \frac{\partial \tilde{f}(z)}{\partial z} (dz + \mu^* d\bar{z} ).
\end{array}
\end{equation}

Hence,
\begin{equation}
    \tilde{f}^*(\rho |dw|^2) = e^{2\lambda_2( \tilde{f}(z) )} |\frac{\partial \tilde{f}(z)}{\partial z}|^2  |dz+\mu^* d\bar{z}|^2.
\end{equation}

So, $\tilde{f}^*(\rho |dw|^2)  = e^{2\lambda_2( \tilde{f}(z) ) - 2\lambda_1( z )} |\frac{\partial \tilde{f}(z)}{\partial z}|^2 g(\mu^*)$. $f^*(\rho |dw|^2)$ is conformal to $g(\mu^*)$. We conclude that $\tilde{f} :(S_1, g(\mu^*))\to (S_2, \rho |dw|^2)$ is conformal. According to Theorem \ref{Ebound}, the Beltrami differential associated to $\tilde{f}$ is the global minimizer of $E_{BC}:\mathcal{T}(S_1,S_2) \to \mathbb{R}$.
\end{proof}

\bigskip

In other words, finding the extremal Teichm\"uller map, $f^*$, is equivalent to minimizing the energy functional $E_{BC}$. During the QC iteration, the Beltrami differential $\mu_n$ is iteratively adjusted and a new map is obtained by $f_n = \mathbf{LBS}(\mu_n)$. It turns out $\mathbf{LBS}(\mu_n)$ is equivalent to computing the harmonic map $f(\mu_n,\rho)$. It can be explained in more details as follows.

\bigskip

\begin{lemma}\label{LBSharmonic}
Suppose $\mu\in \mathcal{T}(S_1,S_2)$. The mapping $f := \mathbf{LBS}(\mu)$ is a harmonic map between $(S_1, g(\mu))$ and  $(S_2, \rho |dw|^2)$.
\end{lemma}
\begin{proof}
Let $\zeta$ be the coordinates of $S_1$ with respect to the metric $g(\mu)$. Let $h$ be the harmonic map between $(S_1, g(\mu))$ and  $(S_2, g(\rho))$. Then $h$ is a critical point of the following harmonic energy:
\begin{equation}\nonumber
E_{harm}(h) = \int_{S_1} \rho(h(\zeta)) (|h_{\zeta}|^2 + |h_{\bar{\zeta}}|^2) dxdy
\end{equation}
Since $f := \mathbf{LBS}(\mu)$, according to the definition, $f$ is the critical point of the following energy functional:
\begin{equation}\nonumber
E_{LBS}(f) = \int_{S_1} \rho(f(z)) (|f_{\bar{z}} - \mu f_z|^2) dxdy
\end{equation}
We will show that the above two energy functionals have the same set of critical points.

Note that $d\zeta = dz + \mu d\bar{z}$, then
\begin{equation}
d\bar{\zeta} = d\bar{z} + \bar{\mu}dz.
\end{equation}

We obtain
\begin{equation}
dz = \frac{1}{1-|\mu|^2}(d\zeta - \mu d\bar{\zeta});\ \ d\bar{z} = \frac{1}{1-|\mu|^2}(-\bar{\mu} d\zeta + d\bar{\zeta}).
\end{equation}

Hence,
\begin{equation}
dz\wedge d\bar{z} = \frac{1}{1-|\mu|^2} d\zeta \wedge d\bar{\zeta};\ \ h_{\bar{\zeta}} = \frac{1}{1-|\mu|^2}(h_{\bar{z}} - \mu h_z).
\end{equation}

Now, the Jacobian $J_h$ of $h$ and the Jacobian $J_f$ of $f$ are given by
\begin{equation}
J_h = |h_{\zeta}|^2 - |h_{\bar{\zeta}}|^2;\ \ J_f = |f_{z}|^2 - |f_{\bar{z}}|^2
\end{equation}

Hence,
\begin{equation}
\begin{split}
E_{harm}(h) & = \int_{S_1} \rho(h(\zeta)) (2 |h_{\bar{\zeta}}|^2 + J_h ) i d\zeta \wedge d\bar{\zeta} \\
& = \int_{S_1} \frac{2}{1-|\mu|^2} \rho(h(z)) |h_{\bar{z}} - \mu h_z|^2 i dz \wedge d\bar{z} + \int_{S_1} \rho(h(\zeta)) J_h  i d\zeta \wedge d\bar{\zeta}
\end{split}
\end{equation}
Since $\mu \in \mathcal{T}(S_1,S_2)$, $|\mu|$ is a constant. Thus,
\begin{equation}
E_{harm}(h) = \frac{2}{1-|\mu|^2} \int_{S_1} \rho(h(z)) |h_{\bar{z}} - \mu h_z|^2 i dz \wedge d\bar{z} + A
\end{equation}
\noindent where $A$ is the surface area of $S_2$. We conclude that $E_{harm}$ and $E_{LBS}$ has the same set of critical points. Since $f$ is a critical point of $E_{LBS}$, $f$ is also a critical point of $E_{harm}$. Hence, $f$ is a harmonic map between $(S_1, g(\mu))$ and  $(S_2, \rho |dw|^2)$.
\end{proof}

The Beltrami differential $\mu_n \in \mathcal{T}(S_1,S_2)$ is iteratively adjusted during the QC iteration. In the next section, we will prove that $E_{BC}(\mu_n)$ monotonically decreases to the global minimizer of $E_{BC}$.

\section{Proof of convergence}\label{sec:convergence}

In this section, we prove the convergence of the Quasi-conformal iteration to the desired Teichm\"uller map.

\bigskip

\begin{lemma}\label{variation}
Suppose $\mu \in \mathcal{B}(S_1,S_2)$ is deformed by
\begin{equation}\nonumber
\mu \to \mu + \epsilon \nu \in \mathcal{B}(S_1,S_2).
\end{equation}
Then, the variation of $E_{BC}$ satisfies:
\begin{equation}\nonumber
E_{BC}(\mu + \epsilon \nu) \leq E_{BC}(\mu) - 4 \mathbf{Re}\int_{S_1} \epsilon\  \Phi(f(\mu,\rho)) \nu \frac{dz_{\mu}\wedge d\bar{z}_{\mu}}{-2i}+ \mathcal{O}(\epsilon^2).
\end{equation}
\noindent where $z_{\mu}$ is the coordinates of $S_1$ under the metric $g(\mu)$.
\end{lemma}
\begin{proof}
Let $\zeta$ be the coordinate of $S_1$ under the metric $g(\mu + \epsilon \nu)$. For simplicity, let $z = z_{\mu}$. Then, we have
\begin{equation}
dz = d{\zeta} - \epsilon \nu d\bar{\zeta};\ \ d\bar{z} = d\bar{\zeta} - \epsilon \bar{\nu} d{\zeta}.
\end{equation}
The area element with respect to $z$ is given by
\begin{equation}
dz\wedge d\bar{z} = d{\zeta}\wedge d\bar{\zeta} - \epsilon \nu d\bar{\zeta}\wedge d\bar{\zeta} - \epsilon \bar{\nu} d\zeta \wedge d\zeta + \epsilon^2 |\nu|^2 d\bar{\zeta}d \zeta.
\end{equation}
Hence,
\begin{equation}
dz\wedge d\bar{z} = d{\zeta}\wedge d\bar{\zeta}+ \epsilon^2 |\nu|^2 d\bar{\zeta}d \zeta.
\end{equation}
Similarly,
\begin{equation}
d\zeta \wedge d\bar{\zeta} = d{z}\wedge d\bar{z}+ \epsilon^2 |\nu|^2 d\bar{z}d z.
\end{equation}

Let $w = f(\mu , \rho)$. Then,
\begin{equation}
\begin{split}
dw & = w_{\zeta} d\zeta + w_{\bar{\zeta}} d\bar{\zeta}\\
 & = w_{z} dz + w_{\bar{z}} d\bar{z}\\
 & = w_{z} (d\zeta - \epsilon \nu d\bar{\zeta})+ w_{\bar{z}}(d\bar{\zeta} - \epsilon \bar{\nu} d\zeta),
\end{split}
\end{equation}
Therefore,
\begin{equation}
\begin{split}
w_{\zeta} \overline{w_{\zeta}} & = (w_z - \epsilon \bar{\nu} w_{\bar{z}})(\overline{w_z} - \epsilon {\nu} \overline{w_{\bar{z}}})\\
 & = |w_z|^2 + \epsilon^2 |\nu|^2 |w_{\bar{z}}|^2 - \epsilon \nu w_z \overline{w_{\bar{z}}} - \epsilon \bar{\nu}\overline{w_z} w_{\bar{z}}.
\end{split}
\end{equation}
Similarly,
\begin{equation}
\begin{split}
w_{\bar{\zeta}} \overline{w_{\bar{\zeta}}} & = (w_{\bar{z}} - \epsilon {\nu} w_{{z}})(\overline{w_{\bar{z}}} - \epsilon \bar{\nu} \overline{w_{{z}}})\\
 & = |w_{\bar{z}}|^2 + \epsilon^2 |\nu|^2 |w_{{z}}|^2 - \epsilon \nu w_{\bar{z}} \overline{w_{{z}}} - \epsilon \bar{\nu}\overline{w_z} w_{\bar{z}}.
\end{split}
\end{equation}

As a result, we get
\begin{equation}
\begin{split}
&E_{BC}(\mu + \epsilon \nu)\leq E_{harm}(w)   = \int_{S_1} \rho(w(\zeta))(|w_{\zeta}|^2 + |w_{\bar{\zeta}}|^2) \frac{d\zeta \wedge d\bar{\zeta}}{-2i} \\
& = \int_{S_1} \rho(w(z))(|w_{z}|^2 + |w_{\bar{z}}|^2) \frac{dz \wedge d\bar{z}}{-2i} - 4\mathbf{Re} \int_{S_1}  \epsilon  \rho(w(z)) w_z \overline{w_{\bar{z}}} \nu \frac{dz \wedge d\bar{z}}{-2i} + \mathcal{O}(\epsilon^2) \\
& = E_{BC}(\mu)  - 4 \mathbf{Re} \int_{S_1} \epsilon\  \rho(w(z)) w_z \overline{w_{\bar{z}}} \nu \frac{dz \wedge d\bar{z}}{-2i} + \mathcal{O}(\epsilon^2)\\
& =  E_{BC}(\mu)  - 4 \mathbf{Re}\int_{S_1} \epsilon\ \Phi(f(\mu,\rho)) \nu \frac{dz_{\mu}\wedge d\bar{z}_{\mu}}{-2i}+ \mathcal{O}(\epsilon^2).
\end{split}
\end{equation}
This completes the proof of the inequality.
\end{proof}

\bigskip

\begin{theorem}\label{thm:direction}
Suppose $\mu\in \mathcal{T}(S_1,S_2)$. For any $\alpha >0$, there exists $w\in \mathcal{B}(S_1,S_2)$ and $\epsilon:S_1\to \mathbb{R}$ such that:

\noindent (i) $\mu + \epsilon w \in \mathcal{T}(S_1,S_2)$;

\noindent (ii) $|\epsilon(p) w(p)| < \alpha$ and $|w(p)| = |\Phi(f(\mu,\rho))(p)|$ for all $p \in S_1$;

\noindent (iii) $\int_{S_1} \epsilon w \Phi(f(\mu,\rho)) \frac{dz_{\mu}\wedge d\bar{z}_{\mu}}{-2i} \geq 0$.
\end{theorem}

\medskip

\begin{proof}
Let $\tilde{k}= |\mu|$ and $\nu = \overline{\Phi(f(\mu,\rho))}$. Pick $\beta\in \mathbb{R}^+$ such that:

\begin{equation}
\beta \sup_{p\in S_1} |\nu(p)|<\alpha /3.
\end{equation}

Consider $\tilde{\mu} = \mu + \beta \nu$.

Suppose:
\begin{equation}
\begin{split}
\Omega_1 &= \{p\in S_1: \arg (\nu) = \arg(\mu)\};\\
\Omega_2 &= \{p\in S_1: \arg (\nu) = -\arg(\mu)\}.
\end{split}
\end{equation}

Let:
\begin{equation}
\gamma = \int_{\Omega_1} |\nu|^2 \frac{dz_{\mu}\wedge d\bar{z}_{\mu}}{-2i} - \int_{\Omega_2} |\nu|^2 \frac{dz_{\mu}\wedge d\bar{z}_{\mu}}{-2i}.
\end{equation}

\smallskip

If $\gamma > 0$, choose $\tilde{k} < k<\sup_{p\in S_1} |\tilde{\mu}(p)|$.

If $\gamma < 0$, choose $\inf_{p\in S_1} |\tilde{\mu}(p)| < k< \tilde{k}$.

If $\gamma = 0$ (including $\Omega_1 = \Omega_2 = \emptyset$), choose $\inf_{p\in S_1} |\tilde{\mu}(p)| < k<\sup_{p\in S_1} |\tilde{\mu}(p)|$.

Let:
\begin{equation}
r = k\frac{\tilde{\mu}}{|\tilde{\mu}|};\ \ w = \frac{r-\mu}{|r-\mu|}|\nu|\ \ \mathrm{and\ }\epsilon = \frac{|r-\mu|}{|\nu|}.
\end{equation}

By definition, $\mu + \epsilon w  = r = k\frac{\tilde{\mu}}{|\tilde{\mu}|} \in \mathcal{T}(S_1,S_2)$. Hence, (i) is satisfied.

Now,
\begin{equation}
|w(p)| = |\nu(p)| = |\Phi(f(\mu,\rho))(p)|\mathrm{\ for\ all\ } p\in S_1.
\end{equation}

Also,
\begin{equation}
\begin{split}
|\epsilon(p) w(p)| &= |r-\mu|\\
&\leq |r-\tilde{\mu}| + |\tilde{\mu} - \mu|\\
&= |r-\tilde{\mu}| + |\beta \nu|\\
& < \frac{2\alpha}{3} + \frac{\alpha}{3} = \alpha.
\end{split}
\end{equation}

Thus, (ii) is also satisfied.

Finally, it is easy to check that:

\begin{equation}
\int_{S_1\backslash (\Omega_1 \cup \Omega_2 )}  \epsilon w \Phi(f(\mu,\rho)) \frac{dz_{\mu}\wedge d\bar{z}_{\mu}}{-2i}\geq 0.
\end{equation}

Now, if $\gamma >0$,
\begin{equation}
\begin{split}
&\int_{\Omega_1}  \epsilon w \Phi(f(\mu,\rho)) \frac{dz_{\mu}\wedge d\bar{z}_{\mu}}{-2i}+\int_{\Omega_2}  \epsilon w \Phi(f(\mu,\rho)) \frac{dz_{\mu}\wedge d\bar{z}_{\mu}}{-2i}\\
&= \int_{\Omega_1}  (k-\tilde{k}) |\nu|^2 \frac{dz_{\mu}\wedge d\bar{z}_{\mu}}{-2i}-\int_{\Omega_2} (k-\tilde{k}) |\nu|^2  \frac{dz_{\mu}\wedge d\bar{z}_{\mu}}{-2i}\\
& = (k-\tilde{k}) \gamma >0.
\end{split}
\end{equation}

If $\gamma <0$,
\begin{equation}
\begin{split}
&\int_{\Omega_1}  \epsilon w \Phi(f(\mu,\rho)) \frac{dz_{\mu}\wedge d\bar{z}_{\mu}}{-2i}+\int_{\Omega_2}  \epsilon w \Phi(f(\mu,\rho)) \frac{dz_{\mu}\wedge d\bar{z}_{\mu}}{-2i}\\
&= - \int_{\Omega_1}  (\tilde{k}-k) |\nu|^2 \frac{dz_{\mu}\wedge d\bar{z}_{\mu}}{-2i}+\int_{\Omega_2} (\tilde{k}-k) |\nu|^2  \frac{dz_{\mu}\wedge d\bar{z}_{\mu}}{-2i}\\
& = -(\tilde{k}-k) \gamma >0.
\end{split}
\end{equation}

We conclude that $\int_{S_1}  \epsilon w \Phi(f(\mu,\rho)) \frac{dz_{\mu}\wedge d\bar{z}_{\mu}}{-2i}\geq 0$ and hence (iii) is satisfied.
\end{proof}

\bigskip

We can now proceed to prove the convergence of the Quasi-conformal iteration.

\bigskip

\begin{theorem}\label{thm:convergence}
Suppose $S_1$ and $S_2$ are open Riemann surfaces with boundaries of the same topology. Given a smooth boundary correspondence $h:\partial S_1\to \partial S_2$, the Quasi-conformal (QC) iteration (\ref{QCiteration}) converges to the unique extremal map, which is also a Teichm\"uller map.
\end{theorem}
\begin{proof}
Suppose the pair $(f_n,\mu_n)$ is obtained at the $n^\mathrm{th}$ iteration. The QC iteration first compute a new quasi-conformal map by $f_{n+1} = \mathbf{LBS}(\mu_{n})$. According to Lemma \ref{LBSharmonic}, $f_{n+1}$ is a harmonic map between $(S_1, g(\mu_{n}))$ and  $(S_2, \rho|dw|^2)$. The Beltrami differential $\nu_{n+1}$ of $f_{n+1}$ can be computed by $\widetilde{\mu}_{n+1} = \mu_n + \beta \mu(f_{n+1},\mu_n)$. $\mu(f_{n+1},\mu_n)$ denotes the Beltrami differential of $f_{n+1}$ under the auxiliary metric with respect to $\mu_n$, namely, $|dz + \mu_n d\bar{z}|^2$. A new Beltrami differential can then be obtained by projecting $\widetilde{\mu}_{n+1}$ onto $\mathcal{T}(S_1,S_2)$ to get
\begin{equation}
\mathcal{P}(\nu_n) =  \mu_n + \epsilon w_n.
\end{equation}
\noindent Here, $w_n:S_1\to \mathbb{C}$ and $\epsilon: S_1\to \mathbb{R}$ is a suitable real function on $S_1$ such that $|\mu_n + \epsilon \mu(f_{n+1},\mu_n)| \equiv k$, where $k$ is a positive constant.

According to Theorem \ref{thm:direction}, by choosing a suitable $k$, we can assume that
\begin{equation}\label{gradientenergy}
\int_{S_1} \epsilon w_n \Phi(f(\mu,\rho)) \frac{dz_{\mu}\wedge d\bar{z}_{\mu}}{-2i} \geq 0.
\end{equation}

$\mathcal{P}(\nu_n)$ is then smoothed out by the Laplacian operator $\mathcal{L}$ with the constraint that it still preserves Equation \ref{gradientenergy}. We get that

\begin{equation}
E_{BC}(\mu_{n+1}) - E_{BC}(\mu_{n}) = - 4 \mathbf{Re}\int_{S_1} \epsilon\Phi(f(\mu_{n},\rho))w_n \frac{dz_{\mu}\wedge d\bar{z}_{\mu}}{-2i}+ \mathcal{O}(\epsilon^2) \leq 0
\end{equation}

\begin{figure*}[t]
\centering
\includegraphics[height=1.75in]{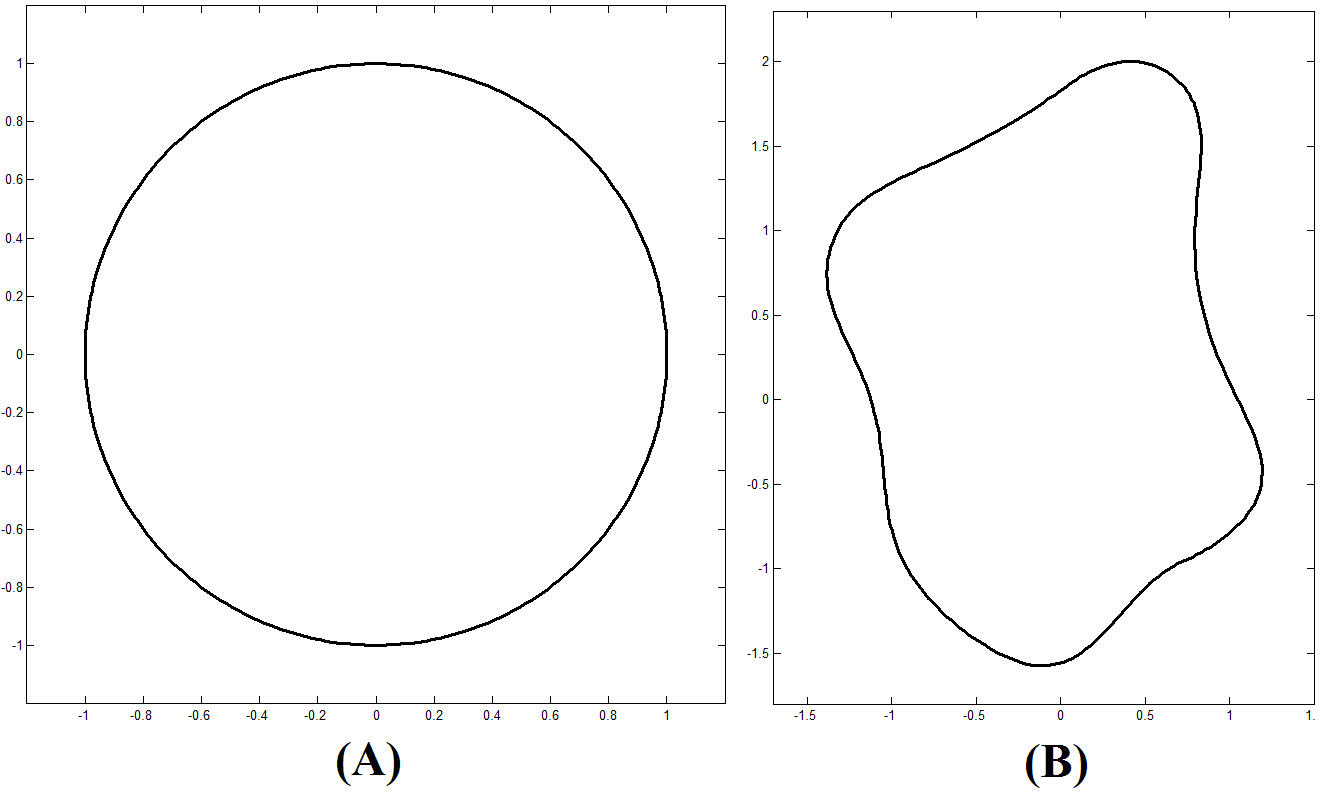}
\caption{Two simply-connected domains. (A) a unit disk $\mathbb{D}$ (B) an arbitrary simply-connected domain. \label{fig:Example1}}
\end{figure*}

\begin{figure*}[t]
\centering
\includegraphics[height=1.75in]{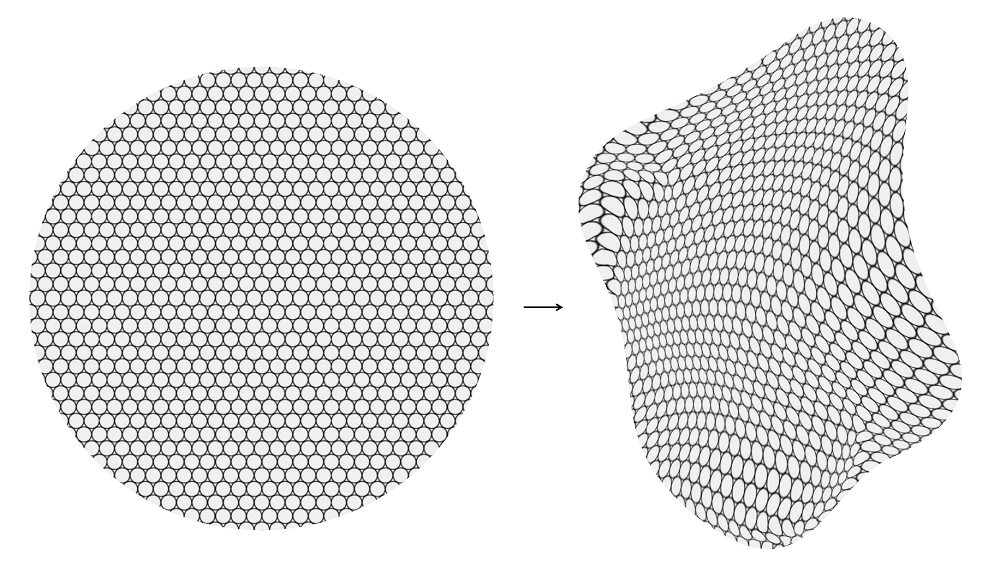}
\caption{Extremal Teichm\"uller map between two simply-connected domains as shown in Figure \ref{fig:Example1}(A) and (B), with given boundary correspondence. \label{fig:Example12}}
\end{figure*}

\begin{figure*}[t]
\centering
\includegraphics[height=1.75in]{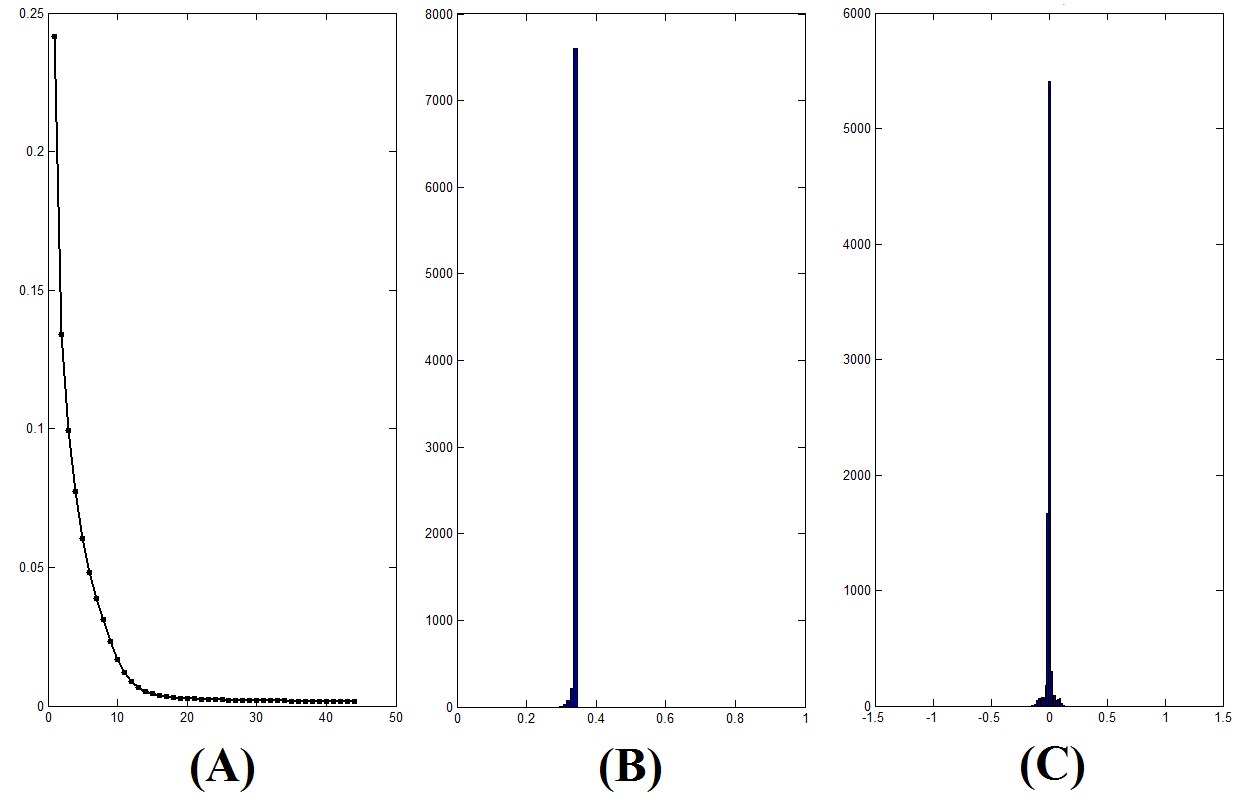}
\caption{(A) shows the energy $E(\mu_n):= E_{BC}(\mu_n) - A(\Omega_2)$ per iterations during the QC iterations of Example 1. (B) shows the histogram of the norm of the optimal Beltrami coefficient $\mu^*$. (C) shows the histogram of the Laplacian of $\mathbf{arg}(\mu^*)$\label{fig:Example13}}
\end{figure*}

Hence, $E(\mu_n)$ is monotonically decreasing. According to Lemma \ref{Ebound}, $E$ is bounded from below. Hence, $E(\mu_n)$ converges. Also, the QC iteration is essential the gradient descend algorithm of $E_{BC}$ and it converges at the critical point $\mu^* = k^* e^{i\theta}$. That is, $\Phi(f(\mu^*,\rho))=0$. In this case, $g(\mu^*)$ is conformal to $\rho$ and hence $f(\mu^*,\rho)$ is a quasi-conformal map with Beltrami differential $\mu^*$.  Furthermore, at the critical point, the Laplacian $\mathcal{L}$ of the Beltrami differential is zero. We conclude that $\theta$ is harmonic. Since $\theta$ is harmonic, we can find its harmonic conjugate $r$ such that $r+ i \theta$ is holomorphic. Define $\varphi = e^{r- i\theta}$, which is also holomorphic. Then, $\mu^* = k^*  \frac{\overline{\varphi}}{|\varphi|}$ is of Teichm\"uller type. Since $\mu^*$ is of Teichm\"uller type, $f(\mu^*,\rho)$ must be a Teichm\"uller map. Now, given a smooth boundary correspondence $h:\partial S_1\to \partial S_2$, there exists a unique Teichm\"uller map which is an extremal map. We conclude that $f(\mu,\rho)$ is the unique extremal Teichm\"uller map.
\end{proof}

\section{Numerical experiments}
Although the numerical testing is not the main focus of this work, we demonstrate some numerical results in this section for the completeness of the paper. The results agree with our theoretical findings.

\subsection*{Example 1} We first test the algorithm to compute the extremal Teichm\"uller map between two simply-connected domains $\Omega_1$ and $\Omega_2$. $\Omega_1$ is chosen to be the unit disk $\mathbb{D}$ as shown in Figure \ref{fig:Example1}(A). $\Omega_1$ is deformed to an arbitrary simply-connected shape $\Omega_2$ as shown in (B). The boundary correspondence $h$ of $\Omega_1$ and $\Omega_2$ is given. We compute the extremal Teichm\"uller map $f:\Omega_1\to \Omega_2$ such that $f|_{\partial \Omega_1} = h$ using the proposed QC iterations. The obtained map is visualized using texture map as shown in Figure \ref{fig:Example12}. The small circles on the source domain is mapped to small ellipses on the target domain with the same eccentricity. Figure \ref{fig:Example13}(A) shows the energy $E(\mu_n):= E_{BC}(\mu_n) - A(\Omega_2)$ versus each iterations in the QC iterations, where $A(\Omega_2)$ is the area of $\Omega_2$. The energy monotonically decreases to 0, which agrees with Theorem \ref{thm:convergence}. (B) shows the histogram of the norm of the optimal Beltrami differential $\mu^*$. It accumulates at 0.33, which illustrates that the obtained map is indeed a Teichm\"uller map. Since $\mu^*$ is of Teichm\"uller type, its argument must be harmonic. (C) shows the histogram of the Laplacian of $\mathbf{arg}(\mu^*)$. It accumulates at 0, meaning that the argument of $\mu^*$ is indeed harmonic.

\subsection*{Example 2} 
\begin{figure*}[t]
\centering
\includegraphics[height=1.75in]{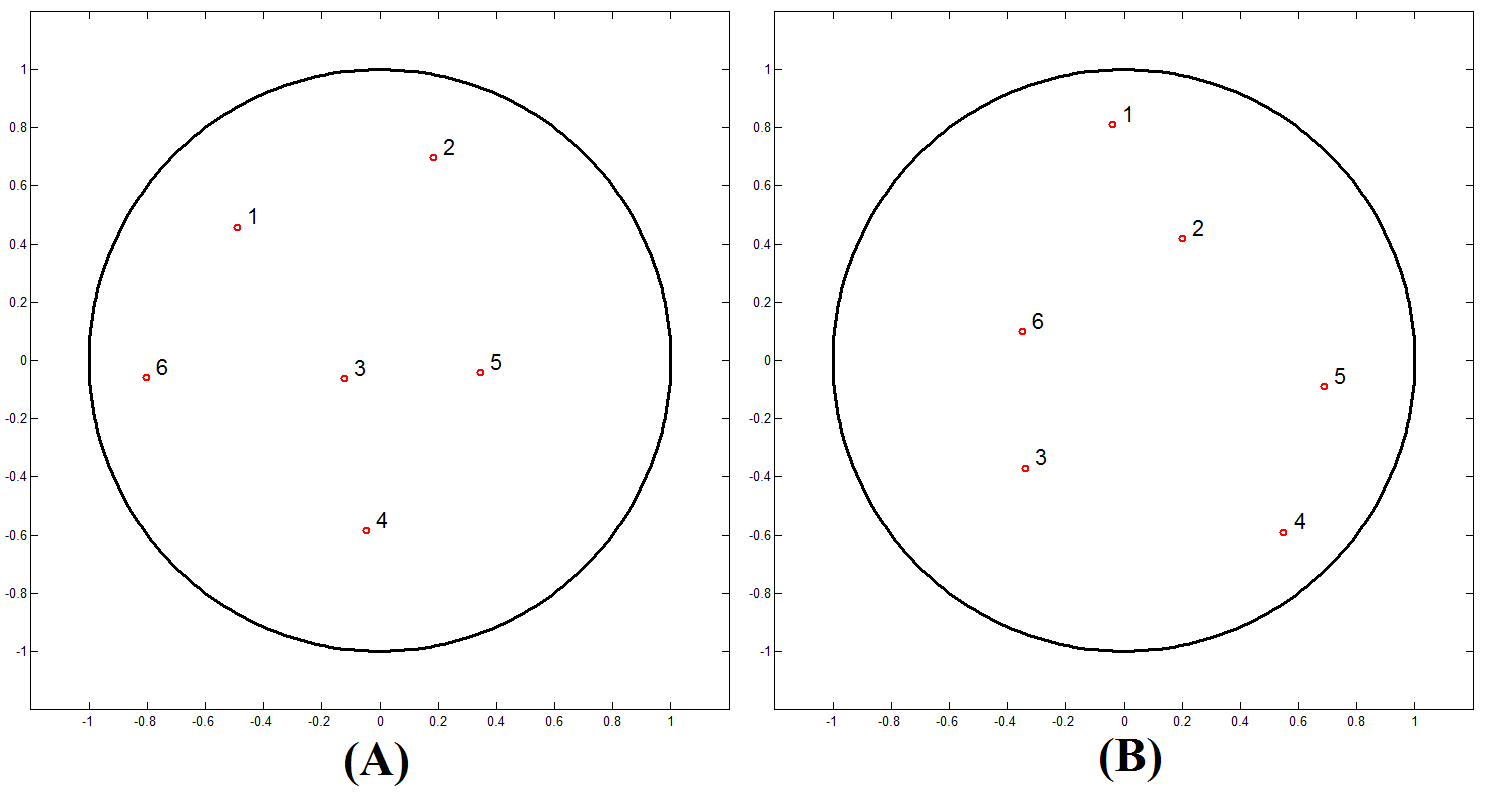}
\caption{Two punctured unit disks. (A) and (B) show two unit disks, each with 6 punctures. \label{fig:Example2}}
\end{figure*}

\begin{figure*}[t]
\centering
\includegraphics[height=1.65in]{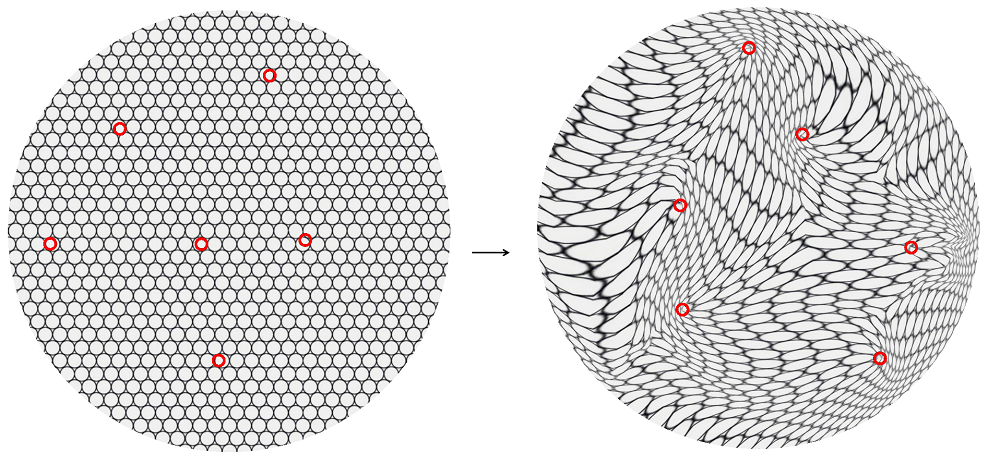}
\caption{Extremal Teichm\"uller map between two punctured unit disks as shown in Figure \ref{fig:Example2}(A) and (B), with given boundary correspondence. . \label{fig:Example22}}
\end{figure*}

\begin{figure*}[t]
\centering
\includegraphics[height=1.75in]{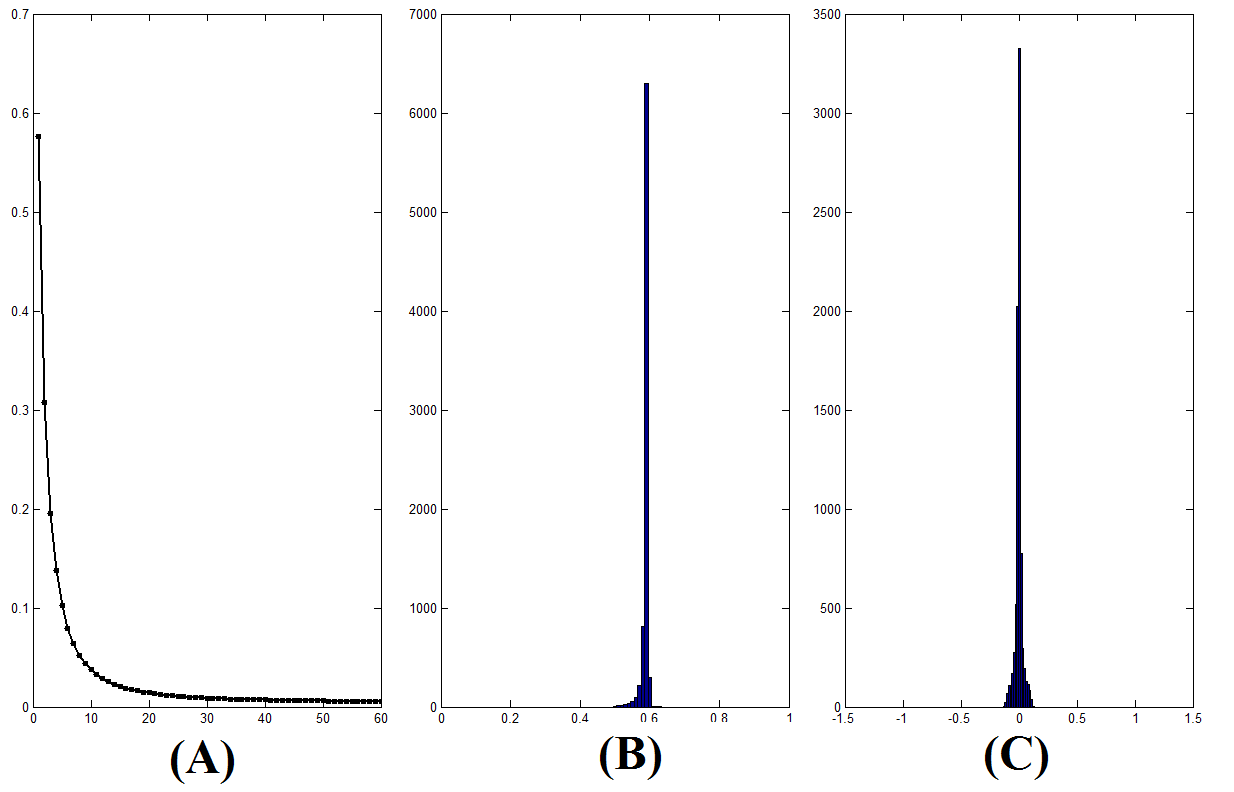}
\caption{(A) shows the energy $E(\mu_n):= E_{BC}(\mu_n) - A(\Omega_2)$ per iterations during the QC iterations of Example 2. (B) shows the histogram of the norm of the optimal Beltrami coefficient $\mu^*$. (C) shows the histogram of the Laplacian of $\mathbf{arg}(\mu^*)$\label{fig:Example23}}
\end{figure*}

In our second example, we test our algorithm to compute the extremal Teichmuller map between two punctured unit disks. Figure \ref{fig:Example2}(A) and (B) show two unit disks, each with 6 punctures. Denote the source domain by $\Omega_1 := \mathbb{D}\setminus \{p_i\}_{i=1}^6$, and denote the target domain by $\Omega_2 := \mathbb{D}\setminus \{q_i\}_{i=1}^6$. The boundary correspondence of $\partial \mathbb{D}$ is chosen to be the identity map. Using the QC iteration, we compute the extremal Teichm\"uller map $f:\Omega_1\to \Omega_2$ such that $f|_{\partial \mathbb{D}} = \mathbf{id}$ and $f(p_i) = q_i$ for $1 \leq i\leq 6$. The obtained map is visualized using texture map as shown in Figure \ref{fig:Example22}. The small circles on the source domain is mapped to small ellipses on the target domain with the same eccentricity. Figure \ref{fig:Example23}(A) shows the energy $E(\mu_n):= E_{BC}(\mu_n) - A(\Omega_2)$ versus each iterations in the QC iterations, where $A(\Omega_2)$ is the area of $Omega_2$. The energy monotonically decreases to 0, which agrees with our theoretical finding. (B) shows the histogram of the norm of the Beltrami differential. It accumulates at 0.6, which illustrates that the obtained map is indeed a Teichm\"uller map. (C) shows the histogram of the Laplacian of $\mathbf{arg}(\mu^*)$. It accumulates at 0, meaning that the argument of $\mu^*$ is indeed harmonic.

\subsection*{Example 3} 
\begin{figure*}[t]
\centering
\includegraphics[height=1.75in]{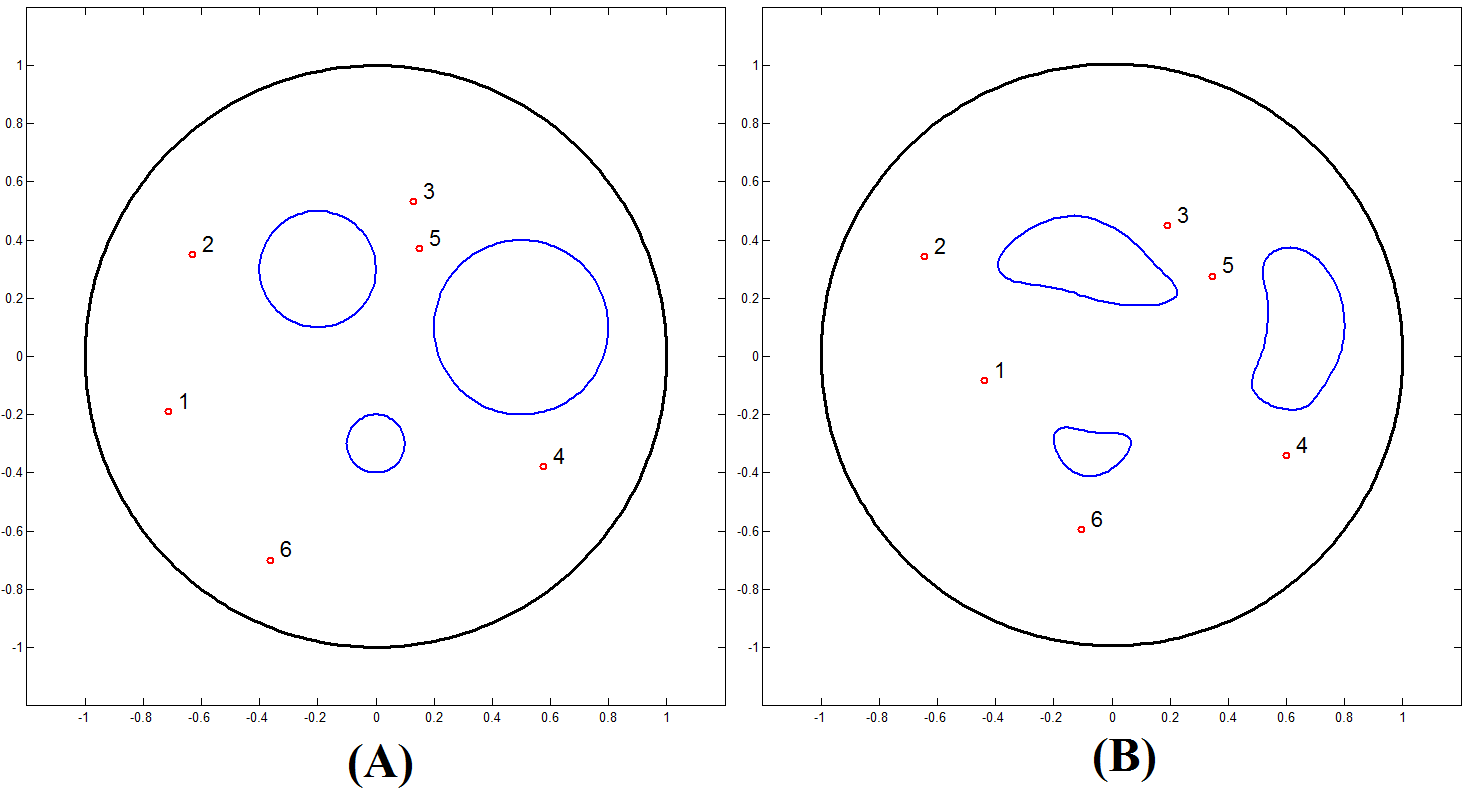}
\caption{Two triply-connected domains, each with 6 punctures. (A) and (B) show the source domain and target domain respectively.\label{fig:Example3}}
\end{figure*}

\begin{figure*}[t]
\centering
\includegraphics[height=1.65in]{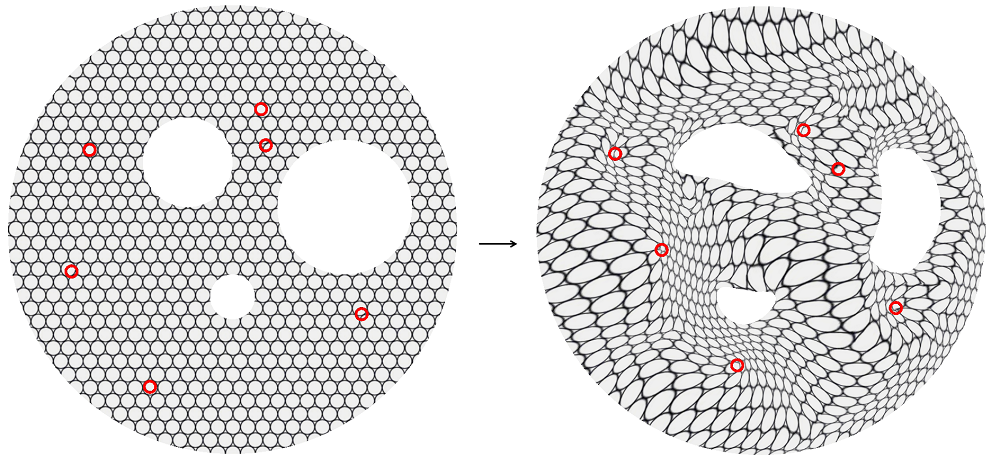}
\caption{Extremal Teichm\"uller map between two triply-connected domains, each with 6 punctures, as shown in Figure \ref{fig:Example3}(A) and (B), with given boundary correspondence. . \label{fig:Example32}}
\end{figure*}

\begin{figure*}[t]
\centering
\includegraphics[height=1.75in]{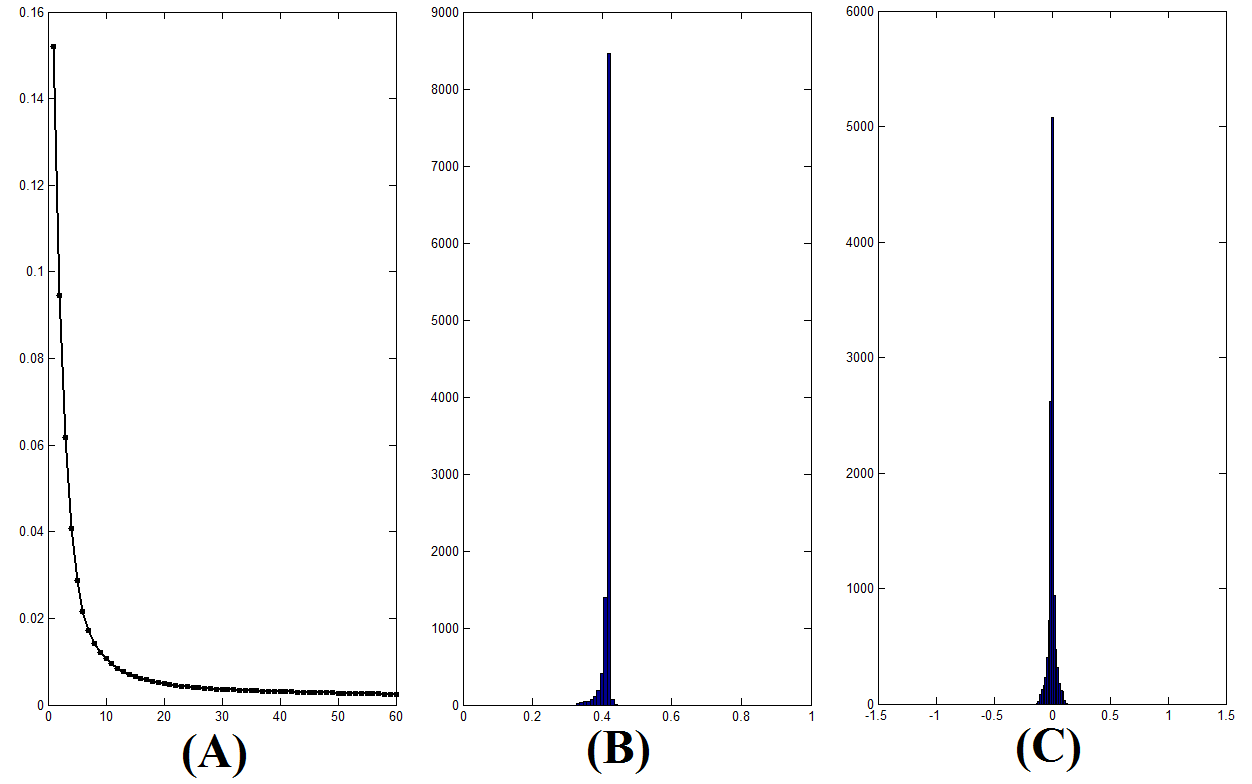}
\caption{(A) shows the energy $E(\mu_n):= E_{BC}(\mu_n) - A(\Omega_2)$ per iterations during the QC iterations of Example 3. (B) shows the histogram of the norm of the optimal Beltrami coefficient $\mu^*$. (C) shows the histogram of the Laplacian of $\mathbf{arg}(\mu^*)$\label{fig:Example33}}
\end{figure*}

Finally, we test the QC iterations to compute the extremal Teichm\"uller map between two triply-connected domains $\Omega_1$ and $\Omega_2$, each with 6 punctures. As shown in Figure \ref{fig:Example3}(A), $\Omega_1$ is chosen to be unit disk with three inner disks and six points removed (denote it by $\{p_i\}_{i=1}^6$). $\Omega_2$ is chosen to be unit disk with three inner regions (with arbitrary shapes) and six points removed (denote it by $\{q_i\}_{i=1}^6$), as shown in (B). Again, the boundary correspondence $h:\partial \Omega_1 \to \partial \Omega_2$ is given. Using the QC iterations, we compute the extremal Teichm\"uller map $f:\Omega_1\to \Omega_2$ such that $f|_{\partial \Omega_1} = h$ and $f(p_i) = q_i$ for $1 \leq i\leq 6$.  The obtained map is visualized using texture map as shown in Figure \ref{fig:Example32}. The small circles on the source domain is mapped to small ellipses on the target domain with the same eccentricity.  Figure \ref{fig:Example33}(A) shows the energy $E(\mu_n):= E_{BC}(\mu_n) - A(\Omega_2)$ versus each iterations in the QC iterations, where $A(\Omega_2)$ is the area of $\Omega_2$. The energy monotonically decreases to 0, which agrees with our theoretical finding. (B) shows the histogram of the norm of the Beltrami differential. It accumulates at 0.42, which illustrates that the obtained map is indeed a Teichm\"uller map. (C) shows the histogram of the Laplacian of $\mathbf{arg}(\mu^*)$. It accumulates at 0, meaning that the argument of $\mu^*$ is indeed harmonic.

\section{Conclusion}\label{sec:conclusion}
This paper gives the convergence proof of the iterative algorithm proposed in \cite{LuiTMap} to compute the extremal Teichm\"uller map between Riemann surfaces of finite type. The iterative algorithm, which is named as quasi-conformal (QC) iteration, can be formulated as the optimization process of the harmonic energy. With this formulation, the QC iteration can be considered as the gradient descent of the harmonic energy under the auxiliary metric given by the Beltrami differentials.

In the future, we will further improve the efficiency of the iterative scheme to optimize the harmonic energy. The proposed framework will also be further extended to compute Teichm\"uller maps between high-genus surfaces (genus $\geq 1$).

\end{document}